\documentclass[11pt,a4paper,abstracton]{scrartcl}
\usepackage[british]{babel}
\usepackage[utf8]{inputenc}
\usepackage{eucal, url}
\usepackage{amsmath,amsxtra, mathrsfs,amsfonts,amssymb}
\numberwithin{equation}{section}

\usepackage{amsthm}

\newtheorem{thm}{Theorem}[section]
\newtheorem{cor}[thm]{Corollary}

\newtheorem{prop}[thm]{Proposition}
\newtheorem{lem}[thm]{Lemma}

\theoremstyle{remark}

\theoremstyle{definition}
\newtheorem{dfn}[thm]{Definition}

\usepackage{paralist}

\def\C{\mathcal{C}}
\def\V{\mathcal{V}}

\def\T{\mathcal{T}}
\def\G{\mathcal{G}}
\def\P{\mathcal{P}}
\def\F{\mathcal{F}}

\def\LO{\mathcal{LO}}

\let\phi\varphi

\DeclareMathOperator{\acl}{acl}
\DeclareMathOperator{\Aut}{Aut}

\newcommand{\ignore}[1]{}
\newcommand{\Nesetril}{Nešetřil}

\newcommand{\Fresse}{Fra\"{i}ss\'{e}}

\title{New Ramsey Classes from Old}

\author{Manuel Bodirsky\thanks{Manuel Bodirsky has received funding from the ERC 
under the European Community's Seventh Framework Programme
(FP7/2007-2013 Grant Agreement no. 257039).
}}

\begin{document}

\maketitle

\begin{abstract}
Let $\C_1$ and $\C_2$ be strong amalgamation classes of
finite structures, with disjoint finite signatures $\sigma$ and $\tau$. Then
$\C_1 \wedge \C_2$ denotes the class of 
all finite $(\sigma \cup \tau)$-structures 
whose $\sigma$-reduct is from $\C_1$ and whose $\tau$-reduct is from $\C_2$. 
We prove that when $\C_1$ and $\C_2$ are Ramsey,
then $\C_1 \wedge \C_2$ is also Ramsey. We also discuss variations of this statement,
and give several examples of new Ramsey classes derived from those general results. 
\end{abstract}

\section{Introduction}
A class of relational structures is a \emph{Ramsey class} if
it satisfies a strong combinatorial property that resembles the statement of Ramsey's theorem.
Surprisingly many classical classes of relational
structures turned out to be Ramsey classes. 
\Nesetril~\cite{RamseyClasses} asked whether one may classify all
Ramsey classes that are closed under induced substructures and have the joint embedding property, and he indicated a link to model-theoretic classification of countably infinite homogeneous structures as an approach to such a classification. 
This program has recently attracted
attention because of a fascinating correspondence between Ramsey classes and the concept of extreme amenability in topological dynamics~\cite{Topo-Dynamics}. We would also like to mention that
Ramsey classes play an important role in
classifications of first-order reducts of homogeneous
relational structures~\cite{BP-reductsRamsey},
and for complexity classification of infinite-domain
constraint satisfaction~\cite{Bodirsky-HDR}.
Establishing
that a class has the Ramsey property is often 
a substantial combinatorial challenge, and we
are therefore interested in general transfer principles that allow to prove the Ramsey property by reducing to known Ramsey classes; this will be the topic of this text.

For structures $A$ and $B$ over the same relational signature, 
let $\binom{B}{A}$ denote the set of
all embeddings of~$A$ into~$B$.
When $f$ is such an embedding, we write $f[A]$ for the copy of $A$ in $B$ that is induced by the image of $A$ under $f$ in $B$.
The partition arrow $C\to(B)^A_r$ means that 
for every function $\chi \colon \binom{C}{A} \rightarrow [r]$ 
(a \emph{colouring}
with $r$ colours) there exists $g\in\binom{C}{B}$ such that $\chi$ is constant on $\binom{g[B]}{A}$. In this case we call $g[B]$ 
a \emph{monochromatic copy} of~$B$ in~$C$.
A class of finite relational structures $\C$ has the \emph{Ramsey property (with respect to embeddings)}\footnote{In some papers, a class $\C$ has the Ramsey property if and only if $\C$ satisfies an analogous property where the partition arrow is not about embeddings, but induced substructures. The two variants are closely related; for a discussion, see~\cite{RamseyClasses}.} if 
for all $A,B \in \C$ and $r \in {\mathbb N}$ there exists a $C \in \C$ such that $C\to(B)^A_r$.
It is easy to see that every class $\C$ with the Ramsey property only contains \emph{rigid} structures,
that is, structures with only one automorphism, the identity. 
Note that an \emph{ordered} structure, that is, a structure that has a strict linear order
as one of its relations, is always rigid. 
A class of relational structures 
that is closed under isomorphisms 
and has the 
Ramsey property is also called a
\emph{Ramsey class}. 

Examples of Ramsey classes 
are 
\begin{itemize}
\item $\LO$, the class of all finite linear orders (this is equivalent to Ramsey's original theorem);
\item the class of all ordered finite graphs (see~\cite{NesetrilRoedlPartite});
\item the class of all ordered $K_n$-free graphs (see~\cite{NesetrilRoedlPartite});
\item the class of all finite partially ordered sets with a linear extension (see~\cite{RamseyClasses};
\item the class of all finite tournaments with an additional linear orde;
\item the class of all finite naturally ordered $C$-relations on a finite set (this is essentially due to~\cite{Mil79}; see~\cite{BodirskyPiguet}).
\end{itemize}
It is of major interest in combinatorics 
to obtain a more systematic understanding of the question which
classes of structures have the Ramsey property. 

\Nesetril\ made the important observation that Ramsey classes that are closed under taking induced substructures are linked with the concept of
\emph{amalgamation} in model theory. 
We say that a class of structures has the \emph{amalgamation property} if for all
$A,B_1,B_2\in \C$ and embeddings $e_1 \colon A \rightarrow B_1$ and $e_2 \colon A \rightarrow B_2$
there exists  
 a $C \in \C$ and embeddings $f_i$ of $B_i$ to $C$ such that $f_1(e_1(a))=f_2(e_2(a))$ for all $a \in A$.  We call $(A,B_1,B_2,e_1,e_2)$ the \emph{amalgamation diagram},
 and $(C,f_1,f_2)$ an \emph{amalgam} of the diagram $(A,B_1,B_2,e_1,e_2)$ (in $\C$).
 If $\C$ has the amalgamation property only for the special case that $A$ is empty,
we say that $\C$ has the \emph{joint embedding property} (here, our assumption that the signature is relational becomes important). 
The mentioned link between Ramsey theory and amalgamation is that 
every class $\C$ of rigid finite relational structures that is closed under isomorphisms and induced substructures, and that has the joint embedding and the Ramsey property also has 
amalgamation property~\cite{RamseyClasses}. Classes of finite structures with countably many non-isomorphic structures that are closed under isomorphisms, induced substructures, and have the amalgamation property are called \emph{amalgamation classes}. 

The \emph{age} of a relational structure $\Gamma$ is the class of all finite structures
that embed into $\Gamma$. 
A structure is \emph{homogeneous} if any isomorphism between finite induced substructures of $\Gamma$ can be extended to an automorphism of $\Gamma$. 
When $\C$ is an amalgamation class, then \Fresse's theorem shows that there exists a countably infinite homogeneous structure $\Gamma$ whose age is $\C$ 
(see e.g.~\cite{Hodges}). 
The structure $\Gamma$ is unique up to isomorphism, and called the \emph{\Fresse-limit} of $\C$;  these homogeneous limit structures
will play an important role in the proof of our main result. 
The significance of \Nesetril's observation is that the transition to countable homogeneous structures brings new tools for the systematic understanding 
of Ramsey classes; and indeed, under some additional assumptions, there are many 
\emph{classification results} for homogeneous structures (such as the classification of all homogeneous directed graphs~\cite{Cherlin}). 

A \emph{strong amalgam} of an amalgamation diagram 
$(A,B_1,B_2,e_1,e_2)$ is an amalgam $(C,f_1,f_2)$ 
such that $f_1(e_1(A)) = f_2(e_2(A)) = f_1(B_1) \cap f_2(B_2)$.
A class $\C$ has \emph{strong amalgamation} if every
amalgamation diagram has a strong amalgam in $\C$.
We say that $\C$ is a \emph{strong amalgamation class} if $\C$ is closed under isomorphisms, induced substructures, and has the strong amalgamation property. An example of a strong amalgamation class is $\LO$. 
Homogeneous structures $\Gamma$ that arise as the \Fresse-limits
of strong amalgamation classes can be characterized via algebraic closure:
in this context, we define 
the algebraic closure $\acl(A)$ of a finite subset $A = \{a_1,\dots,a_n\}$ 
of the domain of $\Gamma$
to be the set of all those elements of $\Gamma$ which lie in finite orbits of the expansion
$(\Gamma,a_1,\dots,a_n)$ of $\Gamma$ by the constants $a_1,\dots,a_n$. 

\begin{prop}[see (2.15) in \cite{Oligo}]\label{prop:acl}
The age of a homogeneous structure $\Gamma$ has strong amalgamation if and only if for any finite subset $A$ of the domain of $\Gamma$, $\acl(A)=A$. 
\end{prop}

\begin{dfn}[from~\cite{Oligo}]
Let $\C_1$ and $\C_2$ be strong amalgamation classes with disjoint 
signatures $\sigma$ and $\tau$. Then
$\C_1 \wedge \C_2$ denotes the class of all finite $(\sigma \cup \tau)$-structures 
whose $\sigma$-reduct is from $\C_1$ and whose $\tau$-reduct is from $\C_2$. 
\end{dfn}

It is clear that $\C_1 \wedge \C_2$ also has strong
amalgamation. In  Section~\ref{sect:full-product}, we prove the following.

\begin{thm}\label{thm:pre-main}
Let $\C_1$ and $\C_2$ be strong amalgamation classes with the Ramsey property, with disjoint finite signatures $\sigma$ and $\tau$.
Then $\C_1 \wedge \C_2$ is Ramsey. 
\end{thm}

The following is an immediate consequence of Theorem~\ref{thm:pre-main}
and the previously known Ramsey results mentioned above. 

\begin{cor}\label{cor:applications}
The following classes of finite structures are Ramsey. 
\begin{enumerate}
\item The class of all permutations of a finite set (represented by two linear orders);
\item The class of all finite sets carrying $n$ linear orders;
\item The class of all finite posets with a linear extension and an additional arbitrary linear order;
\item The class of all finite sets carrying two posets, a linear extension of the first, and a linear extension of the second poset;
\item The class of all finite sets carrying a poset and a linear extension of it, 
and additional linear order and a graph relation.
\item The class of all all naturally ordered $C$-relations on finite sets that additionally carry a poset and a linear extension of this poset. 
\end{enumerate}
\end{cor}

Item 1.\ in Corollary~\ref{cor:applications} has been obtained independently 
by  B\"ottcher and Foniok~\cite{BoettcherFoniok} and by Soki\'c~\cite{Sokic}. 
It is clear that the list can be prolonged easily. 

To prove the statement in item 1., Soki\'c developed a technique called
\emph{cross construction}; also see~\cite{SokicIsrael}. He also proved
item 2.\ and 3. in Corollary~\ref{cor:applications}. The present work has been found independently from~\cite{SokicIsrael}, and it would be interesting to compare our approach with the approach in~\cite{SokicIsrael}.

Homogeneous structures with a finite relational signature
are \emph{$\omega$-categorical}, 
that is, their first-order theory has only one countable model up to isomorphism. 
We can weaken the assumption of having a finite signature slightly, 
and prove the following stronger version which captures 
several additional interesting classes (see Corollary~\ref{cor:examples}). 

\begin{thm}\label{thm:main}
Let $\C_1$ and $\C_2$ be strong amalgamation classes with $\omega$-categorical
\Fresse-limits. If $\C_1$ and $\C_2$ have disjoint relational signatures and the Ramsey property, then $\C_1 \wedge \C_2$ is also Ramsey. 
\end{thm}

To show the Ramsey property for even more classes, we would also like to be able to
generate Ramsey classes that only have one linear order in their signature; this can be accomplished using the following proposition whose proof can be found in Section~\ref{sect:forget-order}.

\begin{prop}\label{prop:forget-order}
Let $\C_1$ and $\LO \wedge \C_2$ be Ramsey classes 
with strong amalgamation and $\omega$-categorical \Fresse-limits,
and suppose that $\C_1$, $\C_2$, and $\LO$ have pairwise disjoint relational signatures.
Then $\C_1 \wedge \C_2$ has the Ramsey property. 
\end{prop}

Proposition~\ref{prop:forget-order} 
is a versatile tool to construct a variety of new Ramsey classes. 
To state many examples, we make the following definitions.
\begin{dfn}
Write 
\begin{itemize}
\item $\T$ for the class of all finite tournaments, 
\item $\G$ for the class of all finite graphs,
\item $\F_n$ for the class of all finite $K_n$-free graphs,
\item $\vec \T$ for the class of all linearly ordered finite tournaments, 
\item $\vec \G$ for the class of all linearly ordered finite graphs,
\item $\vec \F_n$ for the class of all linearly ordered finite $K_n$-free graphs,
\item $\vec \P$ for the class of all linearly extended finite posets,
\item $\vec \C$ for the class of all naturally ordered $C$-relations on a finite set~\cite{BodirskyPiguet},
\item $\vec \V$ for the class of all finite affine vector spaces $V$, equipped with a `natural order' (see~\cite{Topo-Dynamics}); the vector spaces will be
represented as relational structures with an infinite signature that contains a relation symbol for every affine equation.
\end{itemize}
\end{dfn}

\begin{cor}\label{cor:examples}
Let $\C_1$ be one of the classes $\vec\T,\vec\G,\vec\F_n,\vec\P,\vec\C,\vec\V$, and let $\C_2$ be one of the classes $\T,\G,\F_n$.
Then $\C_1 \wedge \C_2$ is Ramsey. 
\end{cor}


\subsection{Topological Dynamics}
Our combinatorial result translates nicely
into a result that shows that certain intersections of extremely amenable groups are again extremely amenable,
based on a connection between Ramsey theory and topological dynamics (Theorem~\ref{thm:kpt}). 
In fact, our presentation of the proof of Theorem~\ref{thm:main} makes use of this connection, and so we briefly present it in the following. 

The property of $\omega$-categoricity of a structure $\Gamma$ can
be characterized in terms of the automorphism group of $\Gamma$.
A countable structure is $\omega$-categorical if and only if
its automorphism group is \emph{oligomorphic}, that is, has only finitely many orbits of $n$-tuples, for all $n$.
A topological group $G$ is called \emph{extremely amenable} if every continuous action of $G$ on a compact Hausdorff space has a fixed point. 
We say that a homogeneous structure $\Gamma$ 
is Ramsey if the class of all finite induced substructures that embed into 
$\Gamma$ is a Ramsey class. 
The following is the central result from~\cite{Topo-Dynamics} (in a slightly 
more general setting --- but we focus on $\omega$-categorical structures here).

\begin{thm}[of~\cite{Topo-Dynamics}]\label{thm:kpt}
Let $\Gamma$ be homogeneous and $\omega$-categorical with domain $D$. Then the following are equivalent.
\begin{itemize}
\item $\Gamma$ is Ramsey.
\item $\Aut(\Gamma)$ is extremely amenable.
\item $\Gamma$ is Ramsey and there is a linear order on $D$ with a quantifier-free first-order definition in $\Gamma$. 
\end{itemize}
\end{thm}
Some explanations are in place, since Theorem~\ref{thm:kpt} is only implicitly in~\cite{Topo-Dynamics}. 
Recall that we color \emph{embeddings}, and not
induced substructures, so the Ramsey property implies rigidity. 
The equivalence of rigidity and having an invariant linear order for Ramsey structures is stated
in Proposition 4.3 in~\cite{Topo-Dynamics}.
Also note that our additional assumption that 
$\Gamma$ is $\omega$-categorical implies that when a relation
(such as the linear order) is preserved
by all automorphisms of $\Gamma$, then it has a first-order definition in $\Gamma$. This
is a well-known consequence of the proof of the theorem of Engeler, Svenonius, and Ryll-Nardzewski; see~\cite{Hodges}. 
Finally, homogeneity of $\Gamma$ implies that first-order formulas are
equivalent to quantifier-free first-order formulas (again, see~\cite{Hodges}).

\section{Model-Complete Cores}
In our proofs, we make use of the concept of \emph{model-complete cores} of $\omega$-categorical structures. 
A structure $\Gamma$ is called 
a \emph{core} if every endomorphism\footnote{An \emph{endomorphism} of $\Gamma$ is a homomorphism from $\Gamma$ to $\Gamma$.} of $\Gamma$ is
an embedding. A first-order theory $T$ is called \emph{model-complete} if all embeddings between models of $T$
preserve all first-order formulas. An $\omega$-categorical structure $\Gamma$ has a model-complete theory if and only if all self-embeddings $e$
of $\Gamma$ are locally generated by the automorphisms of
$\Gamma$, that is, for every finite tuple $t$ of elements from $\Gamma$
there exists an automorphism $\alpha$ of $\Gamma$ such that $e(t) = \alpha(t)$ (see e.g.~Theorem 3.6.11 in~\cite{Bodirsky-HDR}).
In this case, we say that $\Gamma$ is model-complete. 
The following has been shown in~\cite{Cores-Journal} (
also see~\cite{BodHilsMartin-Journal}).

\begin{thm}\label{thm:cores}
Every $\omega$-categorical structure is homomorphically equivalent\footnote{Two structures $\Gamma$ and $\Delta$ are \emph{homomorphically equivalent} if there is a homomorphism from $\Gamma$ to $\Delta$ and a homomorphism from $\Delta$ to $\Gamma$.} to
a model-complete core $\Delta$, which is unique up to isomorphism, and again $\omega$-categorical or finite.
The expansion of $\Delta$ by all existential positive definable relations is homogeneous. 
\end{thm}

The structure $\Delta$ in Theorem~\ref{thm:cores}
will be called \emph{the model-complete core of $\Gamma$}.
We need the following observation.

\begin{prop}\label{prop:mc-core-homogeneous}
The model-complete core $\Delta$ of an $\omega$-categorical homogeneous structure $\Gamma$ is homogeneous.
\end{prop}
\begin{proof}
Let $h$ be a homomorphism from $\Gamma$ to $\Delta$,
and let $i$ be a homomorphism from $\Delta$ to $\Gamma$. 
Suppose that $f$ is an isomorphism between two finite substructures $A,A'$ of $\Delta$.
The restriction of $i$ to $A$ and to $A'$ is an isomorphism as well, since
otherwise the endomorphism $x \mapsto h(i(x))$ of $\Delta$
would not be an embedding, contradicting the assumption that $\Delta$ is a core. 
By homogeneity of $\Gamma$ there exists an automorphism $\alpha$ 
of $\Gamma$ that extends the isomorphism $i \circ f \circ i^{-1}$ between $i(A)$ and $i(A')$. The mapping $e \colon x \mapsto h(\alpha i(x))$ is an endomorphism
of $\Delta$, and therefore an embedding. Since $\Delta$ is a model-complete core, this mapping is locally generated
by the automorphisms of $\Delta$, and in particular there exists an automorphism $\beta$ of $\Delta$ such that 
$\beta(x)=e(x)=f(x)$ for all $x \in A$. This proves homogeneity of $\Delta$.
\end{proof}

We now prove the following.

\begin{thm}\label{thm:mc-core-Ramsey}
Let $\Gamma$ be $\omega$-categorical, homogeneous, and Ramsey, 
and let $\Delta$ be the model-complete core of $\Gamma$. 
Then $\Delta$ is also Ramsey. 
\end{thm}
\begin{proof}
First note that $\Delta$ is homogeneous, by Proposition~\ref{prop:mc-core-homogeneous}.
Let $h$ be a homomorphism from $\Gamma$ to $\Delta$,
and let $i$ be a homomorphism from $\Delta$ to $\Gamma$. 
Write $D$ for the domain of $\Delta$, and $\tau$ for the signature of $\Gamma$ and $\Delta$. 

By Theorem~\ref{thm:kpt}, there is a linear order $<$ on the elements of $\Gamma$ 
with a quantifier-free first-order definition $\phi(x,y)$ in $\Gamma$. We claim that
$\phi$ defines a linear order on the elements of $\Delta$. 
If there were three elements $x,y,z$
such that $\phi(x,y)$, $\phi(y,z)$, and $\phi(z,x)$ in $\Delta$, then $\phi(i(x),i(y))$, $\phi(i(y),i(z))$, and $\phi(i(z),i(x))$ in $\Gamma$ since the restriction of $i$ to $\{x,y,z\}$ is an embedding. Similarly one can verify that the relation defined by $\phi$ in $\Gamma$ is total
and antisymmetric. 

Since $\Delta$ is ordered by a quantifier-free formula, we can prove the statement that 
$\Delta$ is Ramsey by coloring finite induced substructures rather than embeddings.
Let $P,H$ be two finite substructures of $\Delta$. 
Note that $i(P)$ induces in $\Gamma$ a copy of $P$ since otherwise the endomorphism $x \mapsto h(i(x))$ of $\Delta$ would not be an embedding. Moreover, for every copy $Q$ of $P$ in $\Gamma$ we have that $h(Q)$ induces a copy of $P$ in $\Delta$. To see this, let $\alpha$ be the automorphism of $\Gamma$ that maps $i(P)$ to $Q$; such an $\alpha$ exists by homogeneity of $\Gamma$. Then $e \colon x \mapsto h(\alpha i(x))$ must be an embedding, and 
$e(P)=h(Q)$ which proves the claim. 

Let $\chi \colon {\Delta \choose P} \rightarrow \{0,1\}$ be arbitrary. 
We define a map $\xi \colon {\Gamma \choose P} \rightarrow \{0,1\}$ by setting
$\xi(Q) = \chi(h(Q))$ for every copy $Q$ of $P$ in $\Gamma$. 
Since $\Gamma$ is Ramsey, we find a copy $L$ of $H$ in $\Gamma$
such that $\xi$ is constant on ${L \choose P}$. 
By an argument similar as given above, the restriction $h'$ of $h$ to $L$ is an isomorphism,
and the image of $h'$ induces a copy $M$ of $H$ in $\Delta$.
We are left with the task to show that $\chi$ is constant on ${M \choose P}$. 
So let $P_1, P_2$ be two copies of $P$ in $M$. Let $Q_1$, $Q_2$ be the pre-images of
$P_1$ and $P_2$ with respect to the embedding $h'$. Then $Q_1$ and $Q_2$ are copies of $P$ in $L$
and therefore $\chi(Q_1)=\chi(Q_2)$. It follows that $\xi(P_1)=\xi(P_2)$. 
\end{proof}

\section{The Full Product Structure}
\label{sect:full-product}
Let $\Gamma_1$ and $\Gamma_2$ be two structures 
with the same domain $D$ and disjoint signatures $\sigma$ and $\tau$,
respectively. The \emph{full product}  $\Gamma_1 \boxtimes \Gamma_2$
of $\Gamma_1$ and $\Gamma_2$ is a $(\sigma \cup \tau)$-structure
with domain $D^2$ defined as follows. 
For each $k$-ary $R \in \sigma$, the structure $\Gamma_1 \boxtimes \Gamma_2$ has 
the relation 
$$R^{\Gamma_1 \boxtimes \Gamma_2} = \big \{((a_1,b_1),\dots,(a_k,b_k)) \; | \; (a_1,\dots,a_k) \in R^{\Gamma_1}, b_1,\dots,b_k \in D \big \} \; ,$$ 
and for each $k$-ary $R \in \tau$, it has the relation 
$$R^{\Gamma_1 \boxtimes \Gamma_2} = \big \{((a_1,b_1),\dots,(a_k,b_k)) \; | \; (b_1,\dots,b_k) \in R^{\Gamma_2}, a_1,\dots,a_k \in D \big \} \; .$$

Suppose now that $\Gamma_1$ and $\Gamma_2$ are ordered, and
let $G_1,G_2$ be the automorphism group of $\Gamma_1$ and $\Gamma_2$, respectively. Then the product action of the direct product
$G_1 \times G_2$ on $D^2$ equals the automorphism group of 
$\Gamma_1 \boxtimes \Gamma_2$
(using that both $\Gamma_1$ and $\Gamma_2$ are ordered). 

\begin{prop}\label{prop:product-homogeneous} 
Let $\Gamma_1$ and $\Gamma_2$ be ordered homogeneous structures with the same domain $D$ and disjoint signatures.
Then $\Gamma := \Gamma_1 \boxtimes \Gamma_2$ is homogeneous as well. 
\end{prop}
\begin{proof}
Since $\Gamma_1$ and $\Gamma_2$ are ordered, the relation 
$\{((x,y),(u,v)) \; | \; x=u\}$ and the relation $\{((x,y),(u,v)) \; | \; y=v\}$ are 
preserved by isomorphisms between finite substructures of $\Gamma$.
Hence, an isomorphism $\mu$ between finite substructures of $\Gamma$ gives rise to isomorphisms
between finite substructures of $\Gamma_1$ and $\Gamma_2$, respectively. 
Those can be extended to automorphisms $\alpha,\beta$ of $\Gamma_1$ and $\Gamma_2$, by homogeneity.
Then $(x,y) \mapsto (\alpha(x),\beta(y))$ is an automorphism of $\Gamma$ which extends $\mu$.
\end{proof}

The following is also known under the name \emph{product Ramsey theorem}.

\begin{prop}\label{prop:product-Ramsey} 
Let $\Gamma_1$ and $\Gamma_2$ be $\omega$-categorical structures
with the same domain $D$ and disjoint signatures.
When $\Aut(\Gamma_1)$ and $\Aut(\Gamma_2)$ are
extremely amenable, then
the automorphism group of $\Gamma := \Gamma_1 \boxtimes \Gamma_2$ is oligomorphic and extremely amenable. 
\end{prop}
\begin{proof}
It is easy to bound the number of orbits of $n$-tuples 
in $\Gamma$ by the number of orbits of $n$-tuples
of $\Gamma_1$ and $\Gamma_2$, so $\Gamma$ can be seen to be $\omega$-categorical. 
When $G_1$ and $G_2$ are extremely amenable groups, then $G_1 \times G_2$ is extremely amenable as well (see~\cite{Topo-Dynamics}). 
The statement follows
from the observation that $\Aut(\Gamma_1 \boxtimes \Gamma_2)$ is the product action of $\Aut(\Gamma_1) \times \Aut(\Gamma_2)$ on $D^2$
(here we use that $\Gamma_1$ and $\Gamma_2$ are ordered; see~\cite{Bodirsky-HDR} for a more detailed discussion of full products).
\end{proof}

Strong amalgamation will be used via the following lemma.
A relation is called \emph{injective} if it only contains tuples with pairwise distinct entries. 

\begin{lem}\label{lem:homo-strong-amalgamation}
Let $\tau$ be a relational signature, and 
let $\Gamma$ be a homogeneous $\tau$-structure 
such that the class of all finite $\tau$-structures that embed into $\Gamma$ has the strong amalgamation property. Suppose moreover that all relations of
$\Gamma$ are injective. 
Then every finite structure $F$ that homomorphically maps to $\Gamma$
also has an injective homomorphism to $\Gamma$. 
\end{lem}
\begin{proof}
Let $f$ be a homomorphism from $F$ to $\Gamma$ 
such that the range $f(F)$ of $f$ is maximal. 
If $f$ is injective, we are done, otherwise $F$ has elements $u$
and $v$ such that $f(u)=f(v)$. 
Let $A$ be the structure induced by $f(F) \setminus \{f(u)\}$ in $\Gamma$, 
and let $B_1$ and $B_2$ be two disjoint copies of the structure induced 
by $f(F)$ in $\Gamma$. Let $e_1$ be the embedding of $A$ into $B_1$
that maps an element of $f(F) \setminus \{f(u)\}$ to its copy in $B_1$. 
Similarly, there is an embedding $e_2 \colon A \rightarrow B_2$ that maps an element
of $f(F) \setminus \{f(u)\}$ to its copy in $B_2$. By 
strong amalgamation of the age of $\Gamma$, there exist 
embeddings 
$f_1 \colon B_1 \rightarrow \Gamma$ and $f_2 \colon B_2 \rightarrow \Gamma$
such that $f_1[e_1[A]] = f_2[e_2[A]] = f_1[B_1] \cap f_2[B_2]$.
Then the mapping $f' \colon F \rightarrow \Gamma$ defined by $f'(w) = f_1(e_1(f(w)))$
if $w \neq u$, and defined by $f'(w) = f_2(e_2(f(w)))$ 
if $w \neq v$, is well-defined.
To see that it is a homomorphism, note that 
when $R(x_1,\dots,x_n)$ holds in
$F$, then at most one of the $x_i$ can be mapped to
$f(u)$ since the tuples of $R$ in $\Gamma$ have only pairwise
distinct entries. 
Since $f(x) \neq f(y)$ implies that 
$f'(x) \neq f'(y)$, and since moreover $f'(u) \neq f'(v)$,
the function $f'$ also has a larger range than $f$, a contradiction.
\end{proof}

The following is the central lemma connecting the Fraisse-limit of
$\C_1 \wedge \C_2$ with the full product of the Fraisse-limits of $\C_1$ and $\C_2$,
so that we can ultimately use the product Ramsey theorem. 

\begin{lem}\label{lem:link}
Let $\C_1$ and $\C_2$ be strong amalgamation classes of ordered structures, 
with disjoint 
signatures $\sigma$ and $\tau$, such that all relations of $\C_1$ and all
relations of $\C_2$ are injective. 
Let $\Gamma$ be the \Fresse-limit of $\C_1 \wedge \C_2$ with domain $D$,
and suppose that $\Gamma$ is $\omega$-categorical. 
Let $\Gamma_1$ and $\Gamma_2$ be the
$\sigma$- and $\tau$-reduct of $\Gamma$, respectively. 
If $\Gamma_1$ and $\Gamma_2$ are cores, then the following structures are isomorphic.
\begin{enumerate}
\item $\Gamma$ 
\item the substructure induced by $\{(d,d) \; | \; d \in D\}$ in $\Gamma_1 \boxtimes \Gamma_2$
\item the model-complete core of $\Gamma_1 \boxtimes \Gamma_2$
\end{enumerate}
\end{lem}
\begin{proof}
It is straightforward to verify that $d \mapsto (d,d)$ is an isomorphism between $\Gamma$ and the substructure of $\Gamma_1 \boxtimes \Gamma_2$ 
induced by $\{(d,d) \; | \; d \in D\}$. 

To find an isomorphism between $\Gamma$ and the model-complete core of $\Gamma_1 \boxtimes \Gamma_2$,
it suffices to show that $\Gamma$ is a model-complete core, 
and that $\Gamma$ is homomorphically equivalent to $\Gamma_1 \boxtimes \Gamma_2$.
We then use that the model-complete core is unique up to isomorphism (Theorem~\ref{thm:cores}), which gives us the desired isomorphism. 
Model-completeness of $\Gamma$ follows from homogeneity.
To show that $\Gamma$ is a core,
let $e$ be an endomorphism of $\Gamma$. 
Then $e$ is an endomorphism of
the $\sigma$-reduct $\Gamma_1$ of $\Gamma$, and an endomorphism 
of the $\tau$-reduct of $\Gamma_2$ of $\Gamma$.
Since both $\Gamma_1$ and $\Gamma_2$ are cores,
$e$ must be an embedding of $\Gamma$ into
$\Gamma$, which is what we wanted to show.

We finally show that $\Gamma_1 \boxtimes \Gamma_2$ and $\Gamma$ are homomorphically equivalent. For one direction, recall that $\Gamma$ maps to $\Gamma_1 \boxtimes \Gamma_2$ via the mapping $d \mapsto (d,d)$. 
For the other direction, it suffices to show 
that every finite substructure $F$ of $\Gamma_1 \boxtimes \Gamma_2$ homomorphically maps to $\Gamma$, by a standard compactness argument and $\omega$-categoricity of $\Gamma$ (see e.g.~Lemma 3.1.5 in~\cite{Bodirsky-HDR}). 
By Lemma~\ref{lem:homo-strong-amalgamation}, there is 
an \emph{injective} homomorphism $h_1$ from the $\sigma$-reduct of $F$ to $\Gamma_1$ (recall here that the order of $\Gamma_1$ is by assumption strict). 
Similarly, there is an injective homomorphism $h_2$ from the $\tau$-reduct of $F$ into $\Gamma_2$. Let $U$ be the $(\sigma \cup \tau)$-structure with the same domain as $F$, and with
 relations defined as follows: for each $R \in \sigma$ of arity $k$,
 a $k$-tuple $t$ of elements of $U$ is in $R^U$ if and only if $h_1(t)$ is in $R^{\Gamma_1}$. Similarly we define $R^U$ for relations $R \in \tau$,
 with $h_2$ taking the role of $h_1$ and $\Gamma_2$ taking the role of $\Gamma_1$.
Clearly, $h_1$ is an embedding of the $\sigma$-reduct $U_1$ of $U$
 into $\Gamma_1$, and $h_2$ is an embedding of the 
 $\tau$-reduct $U_2$ of $U$ into $\Gamma_2$. Therefore, $U_1 \in \C_1$ and $U_2 \in \C_2$.  By definition of $\C := \C_1 \wedge \C_2$, we have that
$U \in \C$, and there is an embedding $e$ of $U$
into $\Gamma$. Then $e$ is the desired homomorphism from $F$ to $\Gamma$. 
\end{proof}



\begin{proof}[Proof of Theorem~\ref{thm:main}]
Let $\Gamma$ be the \Fresse-limit of $\C_1 \wedge \C_2$,
and let $\Gamma_1$ and $\Gamma_2$ be the $\sigma$- and $\tau$-reduct of $\Gamma$, respectively. Note that $\Gamma_1$ and $\Gamma_2$ are also homogeneous, and $\omega$-categorical by assumption.
Moreover, since $\C_1$ and $\C_2$ have the Ramsey property,
by Theorem~\ref{thm:kpt} there are linear orders $<^{\Gamma_1}$ and
$<^{\Gamma_2}$ with quantifier-free 
first-order definitions $\phi_1$ and $\phi_2$ in $\Gamma_1$ and $\Gamma_2$,
respectively. 

Let $\Gamma_1^*$ be the structure with the same domain as $\Gamma_1$
whose relations are exactly the injective relations that are first-order definable in
$\Gamma_1$. 
Note that this includes in particular the linear order $<^{\Gamma_1}$,
and so $\Gamma_1^*$ is ordered. 
Also note that $\Gamma^*_1$ is homogeneous and a core, and has the same (oligomorphic) automorphism group as $\Gamma_1$, since $\Gamma^*_1$ contains an $n$-ary relation for each orbit of $n$-tuples of distinct elements from $\Gamma$. 
Moreover, observe that algebraic closure only depends
on the automorphism group, and it follows by Proposition~\ref{prop:acl} that
that the age of $\Gamma_1^*$ has strong amalgamation.
We write $\sigma^*$ for the signature of $\Gamma_1^*$.

Analogously we define the structure $\Gamma_2^*$ from $\Gamma_2$;
we choose the signature $\tau^*$ for $\Gamma_2^*$ such that $\tau^*$ is disjoint 
from $\sigma^*$. Finally, let $\Gamma^*$ be the $(\tau^* \cup \sigma^*)$-structure
obtained as the common expansion of $\Gamma_1^*$ and $\Gamma_2^*$. 
Then $\Aut(\Gamma^*) = \Aut(\Gamma)$, because an orbit of an $n$-tuple $t$
in $\Gamma$ is uniquely given by the orbit of $t$ in $\Gamma_1^*$ and
the orbit of $t$ in $\Gamma_2^*$. In particular, 
the automorphism group of $\Gamma^*$ is oligomorphic, and by the theorem
of Ryll-Nardzewski $\Gamma^*$ is $\omega$-categorical. 
Moreover, $\Gamma^*$ is homogeneous and
therefore the \Fresse-limit of its age. 

The groups $\Aut(\Gamma^*_1) = \Aut(\Gamma_1)$ and $\Aut(\Gamma^*_2)=\Aut(\Gamma_2)$ are 
oligomorphic, and extremely amenable by Theorem~\ref{thm:kpt}. 
By Proposition~\ref{prop:product-homogeneous}, $\Gamma^*_1 \boxtimes \Gamma^*_2$ is homogeneous,
and by Proposition~\ref{prop:product-Ramsey}, $\Aut(\Gamma^*_1 \boxtimes \Gamma^*_2)$ is extremely amenable. 
Then Theorem~\ref{thm:mc-core-Ramsey} (again in combination with Theorem~\ref{thm:kpt}) shows that the model-complete core of $\Gamma^*_1 \boxtimes \Gamma^*_2$ has an extremely amenable automorphism group $G$.
By Lemma~\ref{lem:link}, the model-complete core of $\Gamma^*_1 \boxtimes \Gamma^*_2$ is isomorphic
to $\Gamma^*$, and hence $\Aut(\Gamma^*) = \Aut(\Gamma)$ is extremely amenable. 
We conclude by Theorem~\ref{thm:kpt} that $\C_1 \wedge \C_2$ is Ramsey.
\end{proof}

\section{Forgetting one order}
\label{sect:forget-order}
We finally prove Proposition~\ref{prop:forget-order}: let $\C_1$ and $\LO \wedge \C_2$ be Ramsey classes 
with strong amalgamation and $\omega$-categorical \Fresse-limits,
and suppose that $\C_1$, $\C_2$, and $\LO$ have pairwise disjoint relational signatures.
We have to show that $\C_1 \wedge \C_2$ has the Ramsey property. 
\begin{proof}[Proof of Proposition~\ref{prop:forget-order}]
We use the fact that $\C_3 := \C_1 \wedge (\LO \wedge \C_2)$ is Ramsey by Theorem~\ref{thm:main}. Let $\Gamma$ be the \Fresse-limit of $\C_1$.
By Theorem~\ref{thm:kpt}, there is a linear order $<$ on the elements of $\Gamma$ that has a quantifier-free first-order definition $\phi(x,y)$ in $\Gamma$. 

To show that $\C_1 \wedge \C_2$ is Ramsey, let $A$ and $B$ 
be from $\C_1 \wedge \C_2$. Let $A',B'$ be the expansion of $A,B$ by the relation $<$
defined by $\phi$ over $A$ and $B$, respectively. Note that $A',B' \in \C_3$.
Since $\C_3$ is Ramsey, there exists a
$C' \in \C_3$ such that $C' \to(B')^{A'}_r$. 
Let $C$ be the reduct of $C'$ where we drop the relation $<$. 
We claim that $C \to(B)^A_r$. Let $\chi \colon \binom{C}{A} \rightarrow [r]$ be arbitrary. 
We define a coloring $\chi' \colon \binom{C'}{A'} \rightarrow [r]$ as follows. 
Let $e$ be an arbitrary embedding of $A'$ into $C'$. 
Since $A$ is a reduct of $A'$, and $C$ is a reduct of $C'$ with the same signature, 
the mapping $e$ is also an embedding of $A$ into $C$. Therefore, $e$ is in the range 
of $\chi$, and we can define $\chi'(e) := \chi(e)$. 
Since $C' \to(B')^{A'}_r$, there exists an $f \in \binom{C'}{B'}$ such that
$\chi'$ is constant $c$ on $\binom{f[B']}{A'}$. By the same argument as above,
$f$ is also an embedding of $B$ into $C$. 
We claim that $\chi$ is constant on $\binom{f[B]}{A}$. 
Let $e$ be an arbitrary embedding of $A$ into $f[B]$. 
Recall that $A'$ and $B'$ are the expansion of $A$ and $B$ by the relation $<$ defined
by $\phi$. Since embeddings preserve quantifier-free formulas, $e$ preserves in particular $\phi$. Therefore, the mapping $e$ is an embedding of $A'$ into the substructure $f[B']$
of $C'$. In particular, $e$ is in the range of $\chi'$, and $\chi'(e) = c$. 
It follows that $\chi(e) = c$, which concludes the proof that $\chi$ is constant on $\binom{f[B]}{A}$. 
\end{proof}

\paragraph{Acknowledgements.}
Many thanks to Fran\c{c}ois Bossi\`ere, Jan Foniok, Andr\'as Pongr\'acz, and Miodrag Soki\'c for their valuable comments on earlier versions of this text.

\bibliography{../../global}

\end{document}